\documentclass[11pt,twoside]{article}
\usepackage{amsmath, amssymb, mathrsfs, graphicx, bbm}

\def\scr{\mathscr}

\def\az{\alpha}  \def\bz{\beta}
    \def\dz{\delta}
    \def\fz{\varphi}
\def\gz{\gamma}

\def\vz{\varepsilon}

\def\ffz{\Phi}

\def\aq{{\scr A}}

\def\qd{\quad}
\def\qqd{\qquad}

\def\ll{\left}
\def\rr{\right}

\newcommand{\mathsym}[1]{{}}

\def\scr{\mathscr}

\def\leq{\leqslant}
\def\geq{\geqslant}

\font\cms=cmss9 scaled \magstep1

\def\nnd{\noindent}

\def\thm{\nnd\bg{thm1}}
\def\crl{\nnd\bg{crl1}}
\def\lmm{\nnd\bg{lmm1}}
\def\prp{\nnd\bg{prp1}}

\def\xmp{\nnd\bg{xmp1}}

\def\dethm{\end{thm1}}
\def\decrl{\end{crl1}}
\def\delmm{\end{lmm1}}
\def\deprp{\end{prp1}}
\def\dexmp{\end{xmp1}}

\def\prf{\medskip \noindent {\bf Proof}. }
\def\deprf{\quad $\square$ \medskip}
\def\bg{\begin}
\def\be{\bg{equation}}
\def\de{\end{equation}}

\def\dear{\end{eqnarray}}
\def\lb{\label}

\newcommand{\rf}[2]{[\ref{#1}; #2]}

\def\den{\end{enumerate}}

\def\d{\text{\rm d}}

\begin{document}

\allowdisplaybreaks[4]
\thispagestyle{empty}
\renewcommand{\thefootnote}{\fnsymbol{footnote}}

%\noindent {Acta Mathematica Sinica, English Series}\newline
%\noindent {January 2013, Volume 29, Issue 1, pp 1-32 }

\vspace*{.5in}
\begin{center}
{\bf\Large Discrete Hardy-type Inequalities}
\vskip.15in {Zhong-Wei Liao\footnote{Corresponding author. E-mail: zhwliao@mail.bnu.edu.cn (Tel: +008618810561270)}
}
\end{center}
\begin{center} (School of Mathematical Sciences, Beijing Normal University, Beijing 100875, P. R. China)\\
\vskip.1in
\end{center}
\vskip.1in

%\footnotetext{Received May 17, 2012; accepted June 18, 1012}
%\footnotetext{2000 {\it Mathematics Subject Classifications}.\quad 26D10, 60J60, 34L15.}

\bigskip

\noindent {\bf Abstract} \qd This paper studies the Hardy-type inequalities on the discrete intervals. The first result is the variational formulas of the optimal constants. Using these formulas, one may obtain an approximating procedure and the known basic estimates of the optimal constants. The second result, which is the main innovation of this paper, is about the factor of basic upper estimates. An improved factor is presented, which is smaller than the known one and is best possible. Some comparison results are included for comparing the optimal constants on different intervals. 

\vskip.1in


\vskip.1in


\medskip

\section{Introduction}

For given two constants $p$ and $q$ with $1< p\leq q< \infty$, two positive sequences $\mathbf{u}$ and $\mathbf{v}$ on a discrete interval $[1, N]:= \{1, 2, \dots, N \}$ with $N\leq +\infty$, here is the discrete Hardy-type inequality: 
\begin{equation}\lb{Hardy}
\ll[\sum_{n= 1}^N u_n \ll(\sum_{i=1}^n x_i \rr)^q \rr]^{1/q} \leq A \ll(\sum_{n=1}^N v_n x_n^p \rr)^{1/p},
\end{equation}
where $\mathbf{x}$ is an arbitrary non-negative sequence on $[1, N]$. For saving notations, the constant $A$ is assumed to be optimal. 
 
The purpose of this paper is two-fold. First, we give some variational formulas of the optimal constants. The primary applications of the variational formulas are the approximating procedure and the basic estimates. It is necessary to review the advance of the basic estimates in recent research, cf. \cite{Bliss, Chen4, Chen3, Maz'ya, Opic}. In continuous case, the following result, due to B. Opic \rf{Opic}{Theorem 1.14} and V. G. Maz'ya \rf{Maz'ya}{Theorem 1, pp. 42-43}, is well known
\be\lb{basic}
B \leq A \leq \tilde{k}_{q, p} B, 
\de
where $B$ is a quantity described by $N$, $p$, $q$, $\mathbf{u}$ and $\mathbf{v}$, the factor $\tilde{k}_{q, p}$ is a constant defined by $p$ and $q$:
\be \lb{tilde k_qp}
\tilde{k}_{q,p} = \ll(1+ \frac{q}{p^*} \rr)^{1/q} \ll(1+ \frac{p^*}{q} \rr)^{1/p^*}, 
\de
here $p^*$ is the conjugate number of $p$, i.e. $1/p + 1/{p^*}= 1$. In particular, $\tilde{k}_{p,p}= p^{1/p}(p^*)^{1/{p^*}}$. Afterwards, Chen \rf{Chen3}{Theorem 2.1} gets the same conclusions through the variational formulas of the optimal constants. Furthermore, there is an approximating procedure \rf{Chen3}{Theorem 2.2} based on the variational formulas, which can improve the estimates of the optimal constants step by step. In discrete context, when $p=q$, Chen, Wang and Zhang \cite{Chen2} arrive the corresponding variational formulas and basic estimates which are similar to (\ref{basic}), of course, $B$ must be adjusted appropriately in discrete case. When $p\ne q$, Mao \rf{Mao}{Proposition A.1} gets the similar result, but the factor of basic upper estimates is $p^{1/q} (p^*)^{1/p^*}$, which is a little coarser than $\tilde{k}_{q,p}$. Our first destination is to show the corresponding variational formulas in discrete context with the condition of $p\neq q$. Later on, as applications of these formulas, we obtain the basic estimates and the approximating procedure. Overall, these results can be regarded as an extension of the studies in continuous context \cite{Chen3}. 

Second, we study the upper bounds of the basic estimates of the optimal constants in discrete case. Our result is the factor $\tilde{k}_{q, p}$ in (\ref{basic}) can be improved to $k_{q, p}$: 
\be \lb{k_qp}
k_{q,p} = \ll(\frac{r}{B(\frac{1}{r}, \frac{q- 1}{r})} \rr)^{1/p- 1/q}, 
\de
where $B(a, b)= \int_0^1 x^{a-1} (1- x)^{b- 1} \d x$ is the Beta function and $r= q/p -1$. Moreover, our result shows that the factor is best possible and is consistent with the result of continuous case. In continuous case, the improvement has been worked out, cf. \rf{Bennett3}{Theorem 8}, \rf{Manakov}{Theorem 2}, and \rf{Kufner}{pp. 45-47}. The key is the result of Bliss \cite{Bliss}, which gives an integral inequality that the optimal constants can be attained. However, the analogue of the conclusion in the discrete context is nontrivial, as mentioned in \rf{Bennett3}{page 170, two lines above (61)}, \textquotedblleft I have been unable to prove the discrete analogue of Theorem 8\textquotedblright (here the last result is the continuous case). We are lucky to be able to prove this conclusion which constitutes the second part of this paper. When $p=q$, it is well known that the factor $\tilde{k}_{q,p}$ is sharp (see for instance \rf{Hardy_book}{Theorem 326 and 327}). Note that if we allow $q\rightarrow p$, by the identity 
$$
\lim_{r\rightarrow 0^{+}} B^r \ll(\frac{a}{r}, \frac{b}{r} \rr) = \frac{a^a b^b}{(a+ b)^{a+ b}}, \qqd (a,b >0),
$$
we have
\be\lb{k_pp}
k_{p,p}= \tilde{k}_{p,p}= p^{1/p}(p^*)^{1/{p^*}}. 
\de
It means our improved factor is consistent with the original one when $p= q$. Thus, our main results are devoted to the case of $p< q$.

It is a long time for the research about Hardy-type inequalities, which are one of the major themes of harmonic analysis and represent useful tools e.g. in the theory and practice of differential equations, in the theory of approximation etc. From probabilistic consideration, these inequalities are important tools to study the convergence rate of the corresponding processes. These are the origin and motive of this study. 

This paper is organized as follows. The rest of this section, we give some notations and definitions, then illustrate the main results. In Section 2 and Section 3, we prove the conclusions on discrete half line (i.e. $N= \infty$). The case of finite interval (i.e. $N< \infty$) will be handled unitedly in the final section, which gives some comparison results for the optimal constants and their basic estimates on different intervals. 

For the simplicity of illustration, we need some notations. Let $\hat{v}_i = v_i^{1- p^*}, 1\leq i\leq N$.  For any sequence $\mathbf{x}$ on $[1, N]$, define an operator $H$: 
\be\lb{H}
H\mathbf{x}(n)= 
\begin{cases}
0, & n= 0, \\
\sum_{i=1}^n x_i, & n=1, 2, \cdots, N,
\end{cases}
\de
which means the partial summation of $\mathbf{x}$. The following notations are used frequently: 
$$
\az \wedge \bz = \min \{ \az, \bz \}, \qqd \az \vee \bz = \max \{ \az, \bz \}.
$$

Set
$$
\aq[1, N]= \{\mathbf{x}: \text{$x_1>0$, and $x_i\geq 0$ for $i=2, \dots, N$} \}.
$$
In this article, we use the convention $1/0 = \infty$. Then the optimal constant $A$ can be denoted by the following variational formula: 
\begin{equation}\lb{A}
A = \sup_{\mathbf{x} \in \aq[1, N]} \frac{\ll[\sum_{n= 1}^N u_n \ll(\sum_{i=1}^n x_i \rr)^q \rr]^{1/q}}{\ll(\sum_{n= 1}^N v_n x_n^p \rr)^{1/p}} = \sup_{\mathbf{x} \in \aq[1, N]} \frac{\| H\mathbf{x} \|_{l^q(u)}}{\| \mathbf{x} \|_{l^p(v)}}, 
\end{equation}
where $\| \mathbf{x} \|_{l^q(u)}= \ll[ \sum_{n=1}^{N} u_n x_n^q \rr]^{1/q}$ and similarly to $\| \mathbf{x} \|_{l^p(v)}$. 

For upper estimates, define the single summation operator $I^*$ and the double summation operator $I\!I^*$ as:
\begin{align}
I_n^*(\mathbf{x}) &= \frac{\hat{v}_n}{x_n} \ll(\sum_{i=n}^N u_i (H\mathbf{x}(i))^{q/p^*} \rr)^{p^*/q}, \lb{I*} \\
I\!I_n^*(\mathbf{x}) &= \frac{1}{H\mathbf{x}(n)} \sum_{i=1}^n \hat{v}_i \ll(\sum_{j=i}^N u_j (H\mathbf{x}(j))^{q/p^*} \rr)^{p^*/q}, \lb{II*}
\end{align}
with domain $\aq[1, N]$. 

For lower estimates, there are some differences:
\begin{align}
I_n (\mathbf{x})&= \frac{\hat{v}_n}{x_n} \ll(\sum_{i=n}^N u_i (H\mathbf{x}(i))^{q- 1} \rr)^{p^*- 1}, \lb{I} \\
I\!I_n (\mathbf{x})&= \frac{1}{H \mathbf{x} (n)} \sum_{i=1}^n \hat{v}_i \ll(\sum_{j=i}^N u_j (H\mathbf{x}(j))^{q- 1} \rr)^{p^*- 1}. \lb{II} \\
\end{align}
It is easy to see that $I\!I^*= I\!I$ and $I^*= I$ when $p= q$. To avoid the non-summability problem, the domain of $I$ and $I\!I$ have to be modified to:
$$
\aq_0 [1, N]= \ll\{\mathbf{x} \in \aq[1, N]: \sum_{i=1}^N v_i x_i^p < \infty \rr\}.
$$

With these notations, we can give the main conclusions of the variational formulas of the optimal constants.

\thm\lb{Var}
{\cms
The optimal constant $A$ in the Hardy-type inequality {\rm (\ref{Hardy})} satisfies

(i) upper estimates:
\begin{equation}
A\leq \inf_{\mathbf{x}\in \aq[1, N]} \ll(\sup_{n\in [1, N]} I\!I_n^*(\mathbf{x}) \rr)^{1/p^*} = \inf_{\mathbf{x}\in \aq[1, N]} \ll(\sup_{n\in [1, N]} I_n^*(\mathbf{x}) \rr) ^{1/p^*}.
\end{equation}

(ii) lower estimates:
\begin{equation}
\aligned
A&\geq \sup_{\mathbf{x}\in \aq_0 [1, N]} \|\mathbf{x}\|_{l^p(v)}^{p/q -1} \ll(\inf_{n\in [1,N]} I\!I_n (\mathbf{x})\rr)^{{(p-1)}/q} \\
&\geq \sup_{\mathbf{x}\in \aq_0 [1, N]} \|\mathbf{x}\|_{l^p(v)}^{p/q -1} \ll(\inf_{n\in [1,N]} I_n (\mathbf{x})\rr)^{{(p-1)}/q}. \\
\endaligned
\end{equation}
}
\dethm

Using Theorem \ref{Var} on a appropriate test function, we can obtain the basic estimates (\ref{basic}). To be specific, that is:

\crl\lb{basic estimates}
{\cms
The inequality {\rm (\ref{Hardy})} holds for every $\mathbf{x}\in \aq[1, N]$ if and only if
$B< \infty$, where
\be\lb{B}
B= \sup_{n\in [1, N]} \ll( \sum_{i=1}^n \hat{v}_i \rr)^{1/p^*} \ll(\sum_{j=n}^N u_j \rr)^{1/q}.
\de
Moreover, we have
$$
B\leq A\leq \tilde{k}_{q,p}B,
$$
where $\tilde{k}_{q,p}$ is defined as (\ref{tilde k_qp}), which is independent of $\mathbf{u}$, $\mathbf{v}$ and $N$. 
}
\decrl

Roughly speaking, the conclusion of Corollary \ref{basic estimates} is from the first iteration through an appropriate test function. Moreover, we can improve the estimates step by step through multiple iterations on this test function. The following corollary is based on this idea, which is of great significance to numerical computation.

\crl\lb{approximating}
{\cms
(i) Define
$$
\aligned
x^{(1)}_n &= (H \mathbf{\hat{v}} (n))^\az - (H \mathbf{\hat{v}} (n- 1))^\az, \\
x^{(m+ 1)}_n &= \hat{v}_n \ll(\sum_{i= n}^N u_i (H \mathbf{x}^{(m)} (i))^{q/{p^*}} \rr)^{{p^*}/q}, 
\endaligned
$$
where $\az= q/(p^*+ q)$. If $\dz_1:=  \sup_{n\in [1, N]} \ll( I\!I_n (\mathbf{x}^{(1)}) \rr)^{1/p^*}< \infty$, define a sequence as: 
\be
\dz_m= \sup_{n\in [1, N]} \ll( I\!I^*_n(\mathbf{x}^{(m)})\rr)^{1/p^*}, \qd m=2, 3,  \dots. 
\de
Otherwise, define $\dz_m \equiv \infty$. Then $\dz_m$ is a non-increasing sequence (denote $\dz_\infty$ be the limit of $\dz_m$) and we have
$$
A\leq \dz_\infty \leq \cdots \leq \dz_1 \leq \tilde{k}_{q,p} B. 
$$

(ii) Fix $k\in [1, N]$, define
$$
\aligned
y^{(k, 1)}_n &=
\begin{cases}
\hat{v}_n, & 1\leq n\leq k, \\
0, & n> k, 
\end{cases} \\
y^{(k, m+1)}_n &= \hat{v}_n \ll(\sum_{i=n}^N \ll(H \mathbf{y}^{(k, m)} (i) \rr)^{q- 1} \rr)^{p^*- 1},
\endaligned
$$
and define a sequence as
$$
\aligned
\widetilde{\dz}_m&= \sup_{k\in [1, N]} \big\| \mathbf{y}^{(k, m)} \big\|_{l^p (v)}^{p/q -1} \ll(\inf_{n\in [1, N]} I\!I_n \ll( \mathbf{y}^{(k, m)} \rr) \rr)^{{(p- 1)}/q}, \\
\overline{\dz}_m&= \sup_{k\in [1, N]} \frac{\ll[\sum_{n= 1}^N u_n \ll( H \mathbf{y}^{(k, m)} (n) \rr)^q \rr]^{1/q}}{\ll[\sum_{n=1}^N v_n \ll(y_n^{(k, m)} \rr)^p \rr]^{1/p}}. \\
\endaligned
$$
Then we have $A\geq \widetilde{\dz}_m \vee \overline{\dz}_m$ for all $m\geq 1$.
}
\decrl

Another main result of this paper is about the factor in (\ref{basic}). Just like the continuous case, the factor of the basic upper estimates can be improved. Furthermore, we can prove the improved factor is best possible. This result is described in detail below.
\thm\lb{kqpB}
{\cms
The basic upper estimates can be improved to: 
\begin{equation}
A\leq k_{q, p}B, 
\end{equation}
where $k_{q,p}$ is defined as (\ref{k_qp}). In particular, when $N= \infty$ and $\sum_{i=1}^\infty \hat{v}_i = \infty$, the factor $k_{q,p}$ is sharp. 
}
\dethm

\section{Proof of Theorem \ref{Var} }

In this section, we always assume $N= \infty$, and the case of the finite interval will be discussed in Section \ref{Interval}. 

The first proposition is about the property of the sequences which reach the equality case of (\ref{Hardy}). This result is useless in the proof of Theorem \ref{Var}, but is necessary to study the property of the optimal constants. 
 
\prp\lb{decreasing}
{\cms
If the optimal constant $A$, appearing in (\ref{Hardy}), is attained by some non-negative sequence $\mathbf{x}$. Define $w_n = \hat{v}_n^{-1} x_n$, then the sequence $\mathbf{w}$ is decreasing. 
}
\deprp

\prf
With the definition of $\mathbf{w}$, we can rewrite the Hardy- type inequalities (\ref{Hardy}) as:
\be\lb{Hardy*}
\ll[ \sum_{n=1}^\infty u_n \ll(\sum_{i=1}^n \hat{v}_i w_i \rr)^q \rr]^{1/q}\leq A \ll(\sum_{n=1}^\infty \hat{v}_n w_n^p \rr)^{1/p},
\de
and $A$ is attained at $\mathbf{w}$.

The idea in the remainder of this proof is from Bennett \rf{Bennett1}{Section 3}. Using reduction to absurdity, assume there exist integers $i$ and $j$ with $1\leq i< j< \infty$, but $w_i< w_j$. We can construct a new sequence $\mathbf{w'}$ from $\mathbf{w}$ as
\begin{equation}\lb{decreasing1}
w'_i= w'_j= w_0,
\end{equation}
where $w_0$ satisfies 
\begin{equation}\lb{decreasing2}
(\hat{v}_i+ \hat{v}_j) w_0^p= \hat{v}_i w_i^p+ \hat{v}_j w_j^p.
\end{equation}

On the one hand, from (\ref{decreasing2}), the right side of (\ref{Hardy*}) is unchanged when $\mathbf{w}$ is replaced by $\mathbf{w'}$. On the other hand, since $p> 1$ and (\ref{decreasing2}), we have
\be\lb{decreasing3}
w_i< w_0< w_j.
\de
\be\lb{decreasing4}
\hat{v}_i w_i + \hat{v}_j w_j < (\hat{v}_i + \hat{v}_j) w_0.
\de
Combine (\ref{decreasing3}) and (\ref{decreasing4}), we obtain 
$$
\sum_{k=1}^n \hat{v}_k w_k < \sum_{k=1}^n \hat{v}_k w'_k, \qqd \forall n \geq 1. 
$$
It means the left side of (\ref{Hardy*}), increases strictly when $\mathbf{w}$ is replaced by $\mathbf{w'}$. Hence $A$ is not attained at $\mathbf{w}$, which is a contradiction.
\deprf

\noindent {\bf Proof of Theorem \ref{Var}.} The briefing of the proof of Theorem \ref{Var} is given as follows.

(a) First, we need to verify the relation of the single summation operator $I^*$ and the double summation operator $I\!I^*$:
$$
\inf_{\mathbf{x} \in \aq[1, \infty)} \ll[\sup_{n\in [1, \infty)} I\!I_n^*(\mathbf{x}) \rr]^{1/p^*} = \inf_{\mathbf{x} \in \aq[1, \infty)} \ll[\sup_{n\in [1, \infty)} I_n^*(\mathbf{x}) \rr]^{1/p^*}.
$$

For any $\mathbf{x} \in \aq[1, \infty)$, as an application of the proportional property, we get
$$
\aligned
\sup_{n\in [1, \infty)} I\!I_n^*(\mathbf{x}) &= \sup_{n\in [1, \infty)} \frac{1}{H \mathbf{x} (n)} \ll[ \sum_{i= 1}^{n} \hat{v}_i \ll(\sum_{j=i}^{\infty} u_j (H \mathbf{x} (j))^{q/{p^*}} \rr)^{{p^*}/q}\rr] \\
&\leq \sup_{n\in [1, \infty)} \frac{1}{x_n} \ll[ \hat{v}_n \ll(\sum_{i=n}^{\infty} u_i (H \mathbf{x} (i))^{q/{p^*}} \rr)^{{p^*}/q}\rr] \\
&= \sup_{n\in [1, \infty)} I_n^*(\mathbf{x}).
\endaligned
$$
Hence, we have
$$
\inf_{\mathbf{x} \in \aq[1, \infty)} \ll[\sup_{n\in [1, \infty)} I\!I_n^*(\mathbf{x}) \rr]^{1/p^*} \leq\inf_{\mathbf{x} \in \aq[1, \infty)} \ll[\sup_{n\in [1, \infty)} I_n^*(\mathbf{x}) \rr]^{1/p^*}.
$$

On the other hand, for any $\mathbf{x}\in \aq[1, \infty)$, define
$$
y_n= \hat{v}_n \ll(\sum_{i= n}^{\infty} u_i (H \mathbf{x} (i))^{q/{p^*}} \rr)^{{p^*}/q}.
$$
Obviously, $y_n> 0$ on $[1, \infty)$, then $\mathbf{y} \in \aq[1, \infty)$. Again, using the proportional property, we have
$$
\aligned
\sup_{n\in [1, \infty)} I_n^*(\mathbf{y}) &= \sup_{n\in [1, \infty)} \ll[ \frac{\sum_{i=n}^\infty u_i (H \mathbf{y} (i))^{q/{p^*}}}{\sum_{i=n}^\infty u_i (H \mathbf{x} (i))^{q/{p^*}}} \rr]^{{p^*}/q} \\
&\leq \sup_{n\in [1, \infty)} \frac{1}{H \mathbf{x} (n)} \sum_{i= 1}^n \hat{v}_i \ll(\sum_{j= i}^\infty u_j (H \mathbf{x} (j))^{q/{p^*}} \rr)^{{p^*}/q} \\
&= \sup_{n\in [1, \infty)} I\!I_n^*(\mathbf{x}). 
\endaligned
$$
Since $\mathbf{x}$ is arbitrary, we obtain the conclusion we need.

(b) The next step is to show the upper estimates of the optimal constants. Assume $A$ is attained at a non-negative sequence $\mathbf{a}$. For each positive sequence $\mathbf{h}$, as an application of the H\"older inequality and the H\"older-Minkowski inequality, we have
$$
\aligned
\sum_{n= 1}^{\infty} u_n (H \mathbf{a} (n))^q &= \sum_{n= 1}^{\infty} u_n \ll(\sum_{i= 1}^n a_i v_i^{1/p} h_i^{-1} v_i^{-1/p} h_i \rr)^q \\
&\leq \sum_{n= 1}^{\infty} u_n \ll(\sum_{i= 1}^n a_i^p v_i h_i^{-p} \rr)^{q/p} \ll(\sum_{k= 1}^n v_k^{-{p^*}/{p}} h_k^{p^*} \rr)^{q/{p^*}} \\
&\leq \ll\{ \sum_{n= 1}^{\infty} a_n^p v_n h_n^{-p} \ll[\sum_{i= n}^{\infty} u_i \ll(\sum_{k= 1}^i \hat{v}_k h_k^{p^*} \rr)^{q/{p^*}} \rr]^{p/q} \rr\}^{q/p}.
\endaligned
$$
At the last step, we use the H\"older-Minkowski inequality, which needs the condition $p< q$. In particular, when $p= q$, it is Fubini theorem. Now, making a power $1/q$, we get
\begin{align}\lb{==1}
\ll[\sum_{n= 1}^{\infty} u_n (H \mathbf{a} (n))^q \rr]^{1/q} &\leq \ll\{ \sum_{n= 1}^{\infty} a_n^p v_n h_n^{-p} \ll[\sum_{i= n}^{\infty} u_i \ll(\sum_{k= 1}^i \hat{v}_k h_k^{p^*} \rr)^{q/{p^*}} \rr]^{p/q} \rr\}^{1/p} \notag \\
&\leq \sup_{n\in [1, \infty)} \ll[ \frac{1}{h_n^q} \sum_{i=n}^{\infty} u_i \ll(\sum_{k= 1}^i \hat{v}_k h_k^{p^*} \rr)^{q/{p^*}} \rr]^{1/q} \ll(\sum_{j= 1}^{\infty} v_j a_j^p \rr)^{1/p}.
\end{align}
For any $\mathbf{x} \in \aq[1, \infty)$, let 
$$
h_n= \ll( \sum_{i=n}^{\infty} u_i (H \mathbf{x} (i))^{q/{p^*}} \rr)^{1/q}, 
$$
by the proportional property, we have
$$
\aligned
\sup_{n\in [1, \infty)} & \ll[\frac{1}{h_n^q} \sum_{i=n}^{\infty} u_i \ll(\sum_{j= 1}^i \hat{v}_j h_j^{p^*} \rr)^{q/{p^*}} \rr]^{1/q} \\
&\leq \sup_{n\in [1, \infty)} \ll[\frac{1}{H \mathbf{x} (n)} \cdot \sum_{i= 1}^n \hat{v}_i \ll(\sum_{j= i}^{\infty} u_j (H \mathbf{x} (j))^{q/p^*} \rr)^{p^*/q} \rr]^{1/p^*} \\
&= \sup_{n\in [1, \infty)} {I\!I_n^*(\mathbf{x})}^{1/p^*}.
\endaligned
$$
Inserting this formula into (\ref{==1}), we obtain
$$
A = \frac{\ll( \sum_{n= 1}^\infty u_n (H \mathbf{a} (n))^q \rr)^{1/q}}{\ll( \sum_{n= 1}^\infty v_n a_n^p \rr)^{1/p}}\leq \sup_{n\in [1, \infty)} {I\!I_n^*(\mathbf{x})}^{1/p^*}.
$$
Since $\mathbf{x}$ is arbitrary, it follows that
$$
A \leq \inf_{\mathbf{x} \in \aq[1, \infty)} \sup_{n \in [1, \infty)} {I\!I_n^*(\mathbf{x})}^{1/p^*}. 
$$

(c) For the lower estimates, again, we consider the relation of $I$ and $I\!I$ first. For any $\mathbf{x} \in \aq_0 [1, \infty)$, we need to show:
$$
\inf_{n\in [1, \infty)} I_n (\mathbf{x}) \leq \inf_{n\in [1, \infty)} I\!I_n (\mathbf{x}). 
$$
In fact, with the help of proportional property, an argument similar to the one used in (a) can show this result. Next, we should show the variational formulas of $A$.

Since the summability of sequence $\mathbf{x}$, without loss of generality, we may assume $\sum_{i= 1}^\infty v_i x_i^p =1$. Hence our next step is to proof
$$
\sup_{\mathbf{x} \in \tilde{\aq}[1, \infty) } \inf_{n\in [1, \infty)} I\!I_n(\mathbf{x})^{(p- 1)/q} \leq A,
$$
where $\tilde{\aq}[1, \infty)= \{ \mathbf{x} \in \aq_0 [1, \infty): \sum_{i=1}^\infty v_i x_i^p = 1 \}$.

We would begin with the classical variational formulas (\ref{A}) of the optimal constants. For any $\mathbf{x} \in \tilde{\aq}[1, \infty)$, define
$$
y_n= \hat{v}_n \ll(\sum_{i= n}^\infty u_i (H \mathbf{x} (i))^{q-1} \rr)^{p^*- 1},
$$
then we have
\be\lb{==2}
A \geq \frac{\| H \mathbf{y} \|_{l^q (u)}}{\| \mathbf{y} \|_{l^p(v)}} = \frac{\ll[\sum_{n= 1}^\infty u_n (H \mathbf{y} (n))^{q} \rr]^{1/q}}{\ll[\sum_{n= 1}^\infty \hat{v}_n \ll(\sum_{i=n}^\infty u_i (H \mathbf{x} (i))^{q-1} \rr)^{p^*} \rr]^{1/p}}.
\de
Consider the denominator of (\ref{==2}), according to Fubini theorem and the definition of $\mathbf{y}$, we obtain
\begin{align}\lb{==3}
\sum_{i=1}^\infty & y_i \ll(\sum_{j= i}^\infty u_j (H \mathbf{x} (j))^{q-1} \rr) = \sum_{j=1}^\infty u_j (H \mathbf{y} (j)) (H \mathbf{x} (j))^{q-1}\notag \\
&\leq \ll[\sum_{j= 1}^\infty u_j (H \mathbf{y} (j))^{q/p} (H \mathbf{x} (j))^{q/{p^*}} \rr]^{p/q} \ll[\sum_{j= 1}^\infty u_j (H \mathbf{x} (j))^q \rr]^{(q-p)/q}.
\end{align}
The last step is based on the H\"older inequality, which needs the condition $p< q$. Moreover, since $\sum_{i= 1}^\infty v_i x_i^p = 1$ we have
\be\lb{==4}
\ll[\sum_{j= 1}^\infty u_j (H \mathbf{x} (j))^q \rr]^{(q-p)/q} \leq A^{q- p}.
\de
Combining (\ref{==2}), (\ref{==3}), (\ref{==4}) and using the proportional property, we obtain
$$
A \geq \ll[\frac{\sum_{i=1}^\infty u_i (H \mathbf{y} (i))^q}{\sum_{i=1}^\infty u_i (H \mathbf{y} (i))^{q/p} (H \mathbf{x} (i))^{q/{p^*}}} \rr]^{p/{q^2}} \geq \inf_{n \in [1, \infty)} \ll[ \frac{H \mathbf{y} (n)}{H \mathbf{x} (n)} \rr]^{(p-1)/q}.
$$
By the definition of $\mathbf{y}$, we get
$$
\inf_{n\in [1, \infty)} I\!I_n(\mathbf{x})^{(p-1)/q} \leq A.
$$
Since $\mathbf{x}$ is arbitrary, we obtain the variational formulas of $A$. The proof is completed in the case $N= \infty$.
\deprf

\noindent {\bf Proof of Corollary \ref{basic estimates}.} With the help of Theorem \ref{Var}, we can obtain the basic estimates if we choose an appropriate test function.

We consider the upper estimates first. Before the proof, we need some preparations. Given an increasing positive sequence $\mathbf{\ffz}$ on $[1, N]$, for any $n\in [1, N- 1]$ and $0< \az< 1$, we assert
\begin{align}\lb{==5}
\sum_{i=n+ 1}^{N} \ll[\ll(\frac{\ffz_i}{\ffz_n}\rr)^\az- \ll(\frac{\ffz_{i- 1}}{\ffz_n}\rr)^\az \rr] \ll(\frac{\ffz_i}{\ffz_n}\rr)^{-1} \leq \frac{\az}{1- \az}.
\end{align}
In fact, it can be proved by induction. Assume $n=N- 1$. Let $y= \ffz_N/ \ffz_{N- 1}$, then $y\geq 1$ (since $\mathbf{\ffz}$ is increasing). Through simple calculations, we know the function
$$
f(x)= (x^\az- 1) x^{-1}, \qqd x\geq 1,
$$
reaches the maximum when $x= \ll(\frac{1}{1- \az}\rr)^{1/\az}$. Hence
$$
\aligned
\ll[\ll(\frac{\ffz_N}{\ffz_{N- 1}}\rr)^\az- 1 \rr] \ll(\frac{\ffz_N}{\ffz_{N- 1}}\rr)^{-1} &= (y^\az- 1) y^{-1} \\
&\leq (1- \az)^{1/ \az}\ll(\frac{\az}{1- \az}\rr) \\
&\leq \frac{\az}{1- \az}.
\endaligned
$$
For any $1\leq m\leq N- 1$, assume the inequality (\ref{==5}) is true when $n= m$. Now consider $n=m- 1$. Let $y= \ffz_{m}/ \ffz_{m- 1}$, then $y\geq 1$. By the assumption, we have
$$
\aligned
\sum_{i= m}^N & \ll[\ll(\frac{\ffz_i}{\ffz_{m- 1}}\rr)^\az - \ll(\frac{\ffz_{i- 1}}{\ffz_{m- 1}}\rr)^\az \rr] \ll(\frac{\ffz_i}{\ffz_{m- 1}}\rr)^{-1} \\
&= \ll(\frac{\ffz_{m}}{\ffz_{m- 1}}\rr)^{\az- 1} \sum_{i=m}^{N} \ll[\ll(\frac{\ffz_i}{\ffz_{m}}\rr)^\az - \ll(\frac{\ffz_{i- 1}}{\ffz_{m}}\rr)^\az \rr] \ll(\frac{\ffz_i}{\ffz_{m}}\rr)^{-1} \\
&\leq \ll(\frac{\ffz_{m}}{\ffz_{m- 1}}\rr)^{\az- 1} \ll[\frac{\az}{1- \az}+ 1- \ll(\frac{\ffz_{m- 1}}{\ffz_{m}}\rr)^\az \rr] \\
&= \frac{1}{1- \az} y^{\az- 1}- y^{-1}.
\endaligned
$$
Again, by simple calculations, we know the function 
$$
f(x)= \frac{1}{1- \az} x^{\az- 1}- x^{-1}, \qqd x\geq 1,
$$
reachs the maximum $\dfrac{\az}{1- \az}$ when $x= 1$. Hence
$$
\sum_{i= m}^N \ll[\ll(\frac{\ffz_i}{\ffz_{m- 1}}\rr)^\az - \ll(\frac{\ffz_{i- 1}}{\ffz_{m- 1}}\rr)^\az \rr] \ll(\frac{\ffz_i}{\ffz_{m- 1}}\rr)^{-1} \leq \frac{\az}{1- \az}.
$$
By induction, we prove this assertion. We should notice that the right side of (\ref{==5}) is independent of $N$, then (\ref{==5}) also be true when $N \rightarrow \infty$.

Let $\az\in (0, 1)$ be an undetermined parameter. Having the inequality (\ref{==5}) in hand, we can prove:
\be\lb{==6}
\ll(\sum_{i= n}^\infty u_i (H \mathbf{\hat{v}} (i))^{\az q/p^*} \rr)^{1/q}\leq B (H \mathbf{\hat{v}} (n))^{(\az- 1)/p^*} \ll(\frac{1}{1- \az}\rr)^{1/q}, \qd 1\leq n< \infty.
\de
With the definition of $B$, for any $n$, we have $(\sum_{i=n}^\infty u_i)^{1/q}\leq B (H \mathbf{\hat{v}} (n))^{-1/{p^*}}$. Let $\ffz_n= (H \mathbf{\hat{v}} (n))^{q/{p^*}}$, summation by parts, we have
$$\aligned
\sum_{i= n}^\infty u_i \ffz_i^\az &= \ffz_n^\az \ll(\sum_{i=n}^\infty u_i \rr) + \sum_{i=n+ 1}^\infty \ll(\ffz_i^\az- \ffz_{i- 1}^\az \rr) \ll(\sum_{j=i}^\infty u_j \rr) \\
&\leq B^q \ffz_n^{\az- 1} + B^q \sum_{i=n+ 1}^\infty \ll(\ffz_i^\az- \ffz_{i- 1}^\az \rr) \ffz_i^{-1} \\
&= B^q \ffz_n^{\az- 1} \ll\{ 1+ \sum_{i=n+ 1}^\infty \ll[\ll(\frac{\ffz_i}{\ffz_n}\rr)^\az- \ll(\frac{\ffz_{i- 1}}{\ffz_n}\rr)^\az \rr] \ll(\frac{\ffz_i}{\ffz_n}\rr)^{-1} \rr\} \\
&\leq B^q \ffz_n^{\az- 1} \ll(\frac{1}{1- \az} \rr).
\endaligned$$
Making a power $1/q$, we obtain the inequality (\ref{==6}).

Now, we are ready to prove the upper bounds of the basic estimates. For any $n\in[1, \infty)$, Let 
$x_n = (H \mathbf{\hat{v}} (n))^{\az} - (H \mathbf{\hat{v}} (n- 1))^{\az}$, we have
$$
\aligned
{I_n^*(\mathbf{x})}^{1/p^*} &= \ll[\frac{\hat{v}_n}{(H \mathbf{\hat{v}} (n))^{\az} - (H \mathbf{\hat{v}} (n- 1))^{\az}} \ll(\sum_{i= n}^\infty u_i (H \mathbf{\hat{v}} (i))^{\az q /p^*} \rr)^{{p^*}/q} \rr]^{1/{p^*}} \\
&\leq \ll[\frac{1}{\az (H \mathbf{\hat{v}} (n))^{\az- 1} } \ll(\sum_{i= n}^\infty u_i (H \mathbf{\hat{v}} (i))^{\az q/ p^*} \rr)^{{p^*}/q} \rr]^{1/{p^*}} \\
&= \az^{-1/{p^*}} \ll(H \mathbf{\hat{v}} (n) \rr)^{(1- \az)/{p^*}} \ll(\sum_{i= n}^\infty u_i (H \mathbf{\hat{v}} (i))^{\az q/{p^*}} \rr)^{1/q} \\
&\leq B \az^{-1/{p^*}} (1- \az)^{-1/q}.
\endaligned
$$
An easy calculation shows that the function
$$
f(x)= x^{-1/{p^*}} (1- x)^{-1/q}, \qqd 0< x<1, 
$$
reaches the maximum
$$
\tilde{k}_{q, p}= \ll(1+ \frac{q}{p^*}\rr)^{1/q} \ll(1+ \frac{p^*}{q}\rr)^{1/p^*}
$$
when $x= \dfrac{q}{p^*+ q}$. Hence we take $\az= \dfrac{q}{p^*+ q}$, then we have
$$
\sup_{n\in [1, \infty)} \ll( I_n^* (\mathbf{x}) \rr)^{1/p^*}\leq \tilde{k}_{q, p} B.
$$
By Theorem \ref{Var}, we get the basic upper estimates.

The basic lower estimates are more straightforward. For any $n\in [1, \infty)$, we can choose a test sequence as 
$$
x_i^{(n)}= 
\begin{cases}
\hat{v}_i, & 1\leq i\leq n, \\
0, & n< i< \infty. 
\end{cases}
$$
It is obvious that $\mathbf{x}^{(n)} \in \aq_0 [1, \infty)$, then by (\ref{A}) we have
$$
\aligned
A&\geq \sup_{n\in [1, \infty)} \frac{\ll[\sum_{i= 1}^\infty u_i (H \mathbf{x}^{(n)} (i))^q \rr]^{1/q}}{\ll[ \sum_{i= 1}^\infty v_i \ll(x_i^{(n)} \rr)^p \rr]^{1/p}} \\
&= \sup_{n\in [1, \infty)} \ll(\sum_{i= 1}^n \hat{v}_i \rr)^{1/p^*} \ll[ \ll(\sum_{i= 1}^n \hat{v}_i \rr)^{-q} \ll( \sum_{i= 1}^{n- 1} u_i \ll(H \mathbf{\hat{v}} (i) \rr)^q \rr) + \sum_{i= n}^\infty u_i \rr]^{1/q} \\
&\geq B.
\endaligned
$$
This completes the proof of Corollary \ref{basic estimates}.
\deprf

\medskip
\noindent {\bf Proof of Corollary \ref{approximating}.} By the proportional property, we can obtain the monotonicity of \{$\dz_n$\}. The approximating sequence $\{ \dz_n \}$ comes from the upper estimates of the variational formula, and $\{\widetilde{\dz}_n\}$ comes from the lower one. These results are the simple applications of Theorem \ref{Var}. The sequence $\{\overline{\dz}_n\}$ is the straightforward application of the classical variational formula (\ref{A}).
\deprf

\section{Proof of Theorem \ref{kqpB}}

Again, we assume $N= \infty$, and the case of finite interval will be discussed in Section \ref{Interval}. We begin with the following well-known lemma. 

\lmm\lb{compare}
{\cms
Let $\mathbf{a}$, $\mathbf{b}$ be sequences with non-negative entries. If
$$
\sum_{k=i}^\infty a_k \leq \sum_{k=i}^\infty b_k, \qqd (\forall ~ i=1,2, \cdots)
$$
then for any increasing non-negative sequence $\mathbf{c}$, we have
$$
\sum_{k=1}^\infty a_k c_k \leq \sum_{k=1}^\infty b_k c_k.
$$
}
\delmm

\prf
Set $c_0= 0$. Summation by parts, we have
$$
\aligned
\sum_{k=1}^\infty a_k c_k &= \sum_{k=1}^\infty \ll(\sum_{i=1}^k (c_i- c_{i-1}) \rr) a_k = \sum_{i=1}^\infty \ll(\sum_{k=i}^\infty a_k \rr) (c_i- c_{i-1})  \\
&\leq \sum_{i=1}^\infty \ll(\sum_{k=i}^\infty b_k \rr) (c_i- c_{i-1})= \sum_{k=1}^\infty b_k c_k.
\endaligned
$$ 
This completes the proof of Lemma \ref{compare}. 
\deprf

The conclusion of Lemma \ref{compare} is about increasing sequences, analogously, there is a conclusion corresponding to the decreasing sequences, cf. \rf{Bennett2}{Lemma 1}. The following lemma is due to Bliss \cite{Bliss}, which gives a special Hardy-type inequality in the continuous case.
\lmm\lb{Bliss_lemma}
{\cms
For any non-negative real function $f(x)$, we have
\be
\ll(\int_0^\infty \frac{1}{x^{q- r}} \ll(\int_0^x f(t) \d t \rr)^q \d x \rr)^{1/q}\leq k_{q,p} \ll(\frac{p^*}{q} \rr)^{1/q} \ll(\int_0^\infty f^p(x) \d x \rr)^{1/p}, 
\de
where $r= q/p -1$ and $k_{q,p}$ is the optimal constant, which is defined as (\ref{k_qp}). Moreover, the optimal constant is attained when 
$$
f(x)= \frac{c}{(d \cdot x^r + 1)^{(r+ 1)/r}}, 
$$
where $c$ and $d$ are non-negative constants. 
}
\delmm

\medskip
\noindent { {\bf Proof of Theorem \ref{kqpB}}.  By Corollary \ref{basic estimates}, it is obvious that $A= \infty$ if $B= \infty$. To avoid this trivial case, we assume $B< \infty$. 

(a) First we consider the case that $H \mathbf{\hat{v}} (\infty)= \lim_{n \rightarrow \infty} H \mathbf{\hat{v}} (n) =\infty$. Similar to Proposition \ref{decreasing}, we can rewrite the Hardy-type inequalities (\ref{Hardy}) as:
$$
\sum_{n=1}^\infty u_n \ll(\sum_{i=1}^n \hat{v}_i x_i \rr)^q \leq A^q \ll(\sum_{n=1}^\infty \hat{v}_n x_n^p \rr)^{q/p}.
$$

Define sequence $\mathbf{\tilde{u}}$
\be\lb{tilde_u}
\tilde{u}_n = B^q \ll( (H \mathbf{\hat{v}} (n))^{-q/p^*}- (H \mathbf{\hat{v}} (n+1))^{-q/p^*} \rr), \qqd n\geq 1. 
\de
By direct summation and $H \mathbf{\hat{v}} (\infty)= \infty$, we have 
$$
\sum_{i= n}^\infty \tilde{u}_i = B^q \ll(H \mathbf{\hat{v}} (n)^{-q/p^*} - H \mathbf{\hat{v}} (\infty)^{-q/p^*} \rr)= B^q H \mathbf{\hat{v}} (n)^{-q/p^*} \geq \sum_{i= n}^\infty u_i. 
$$
Applying Lemma \ref{compare}, for any non-negative sequence $\mathbf{x}$, we obtain
\be\lb{kqp1}
\sum_{n= 1}^\infty u_n \ll(\sum_{i=1}^n \hat{v}_i x_i \rr)^q \leq \sum_{n= 1}^\infty \tilde{u}_n \ll(\sum_{i=1}^n \hat{v}_i x_i \rr)^q. 
\de

The next is to show
$$
\sum_{n= 1}^\infty \tilde{u}_n \ll(\sum_{i=1}^n \hat{v}_i x_i \rr)^q \leq A^q \ll(\sum_{n=1}^\infty  \hat{v}_n x_n^p \rr)^{q/p}. 
$$
In order to use Lemma \ref{Bliss_lemma}, we should construct a function which connects summation with integration. Defined function $f: [0, \infty) \rightarrow [0, \infty)$
\be\lb{prp1}
f(x)= 
\begin{cases}
x_n, & H \mathbf{\hat{v}} (n-1) \leq x< H \mathbf{\hat{v}} (n), \\
0, & x\geq \sup_n H \mathbf{\hat{v}} (n).
\end{cases}
\de
It is clear that 
\be\lb{prp2}
\sum_{i=1}^\infty \hat{v}_i x_i^p = \int_{0}^\infty f^p(x) \d x, 
\de
and
\be\lb{prp3}
\sum_{i=1}^n \hat{v}_i x_i\leq \int_0^\az f(t) \d t,
\de
where $H \mathbf{\hat{v}} (n) \leq \az< H \mathbf{\hat{v}} (n+ 1)$. 

For convenience, write $\tilde{u}_0= 0, \hat{v}_0= 0$. Applying (\ref{prp3}), Lemma \ref{Bliss_lemma} and (\ref{prp2}), we see that
\begin{align}
\allowdisplaybreaks
\sum_{n=0}^\infty \tilde{u}_n \ll(\sum_{k=0}^n \hat{v}_k x_k \rr)^q &= \sum_{n=0}^\infty B^q \ll( (H \mathbf{\hat{v}} (n))^{-q/p^*}- (H \mathbf{\hat{v}} (n+ 1))^{-q/p^*} \rr) \ll(\sum_{k= 0}^n \hat{v}_k x_k \rr)^q \notag \\
&= \sum_{n=0}^\infty B^q \ll(\frac{q}{p^*} \int_{H \mathbf{\hat{v}} (n)}^{H \mathbf{\hat{v}} (n+ 1)} x^{-q/p^* - 1} \d x\rr) \ll(\sum_{k= 0}^n \hat{v}_k x_k \rr)^q \notag \\
&\leq \frac{q}{p^*} B^q \ll(\sum_{n=0}^\infty \int_{H \mathbf{\hat{v}} (n)}^{H \mathbf{\hat{v}} (n+ 1)} x^{-q/p^* - 1} \ll(\int_0^x f(t) \d t\rr)^q \d x \rr) \notag \\
&= \frac{q}{p^*} B^q \int_0^\infty x^{r- q} \ll(\int_0^x f(t) \d t \rr)^q \d x \qd (r=q/p - 1)\notag \\
&\leq B^q k_{q, p}^q \ll(\int_0^\infty f^p(x) \d x \rr)^{q/p} \notag \\
&= B^q k_{q, p}^q \ll(\sum_{i=0}^\infty \hat{v}_i x_i^p \rr)^{q/p}. \notag
\end{align}
By the definition of the optimal constant, we have $A\leq k_{q, p} B$.

(b) To show the factor of the basic upper estimate is best possible, we attempt to mimic the extremal function in Lemma \ref{Bliss_lemma}: 
$$
f(x)= \frac{c}{\ll(dx^r+ 1\rr)^{(r+1)/r}} \qd \text{and} \qd \int_0^x f(t) \d t = \frac{cx}{\ll(dx^r +1 \rr)^{1/r}},
$$
where $c$ and $d$ are arbitrary positive constants. We start with this form, set
$$
u_n= n^{-q/p^*}- (n+1)^{-q/p^*}, \qqd v_n \equiv 1,
$$
and 
$$
x_n = \frac{c n}{(n^r+ d)^{1/r}}- \frac{c (n-1)}{((n-1)^r+ d)^{1/r}}. 
$$
Obviously, the form of $\mathbf{x}$ comes from the difference of $\int_0^x f(t) \d t$. In this case, we have $B= 1$. Here we are free to choose $c$ and $d$, however, no matter what choices, there is some loss of precision between integrals and series. But by direct calculation, we find that this loss becomes negligible when $c/d \rightarrow 0$. Without loss of generality, we choose $c= 1$ and let $d$ be a positive and large enough real number. Next, we calculate the left and the right side of (\ref{Hardy}). 

The calculation of the right side of (\ref{Hardy}) is direct. By the definition of $\mathbf{x}$, we have
\begin{align}\lb{right}
\sum_{n= 1}^\infty x_n^p &= \sum_{n=1}^\infty \ll[\frac{n}{(n^r + d)^{1/r}} -\frac{n- 1}{((n- 1)^r + d)^{1/r}} \rr]^p \notag \\
&= \sum_{n=1}^\infty \ll[\int_{n- 1}^n \frac{d}{(x^r+ d)^{1/r +1}} \d x \rr]^p \notag \\
&\leq \sum_{n= 1}^\infty \int_{n- 1}^n \ll(\frac{d}{(x^r+ d)^{1/r +1}} \rr)^p \d x \notag \\
&= r^{-1} d^{(1- p)/r} B\ll(\frac{1}{r}, \frac{q- 1}{r}\rr).
\end{align}

The left side of (\ref{Hardy}) is difficult. First, we assert there is a large enough integer $N$ such that 
\be\lb{left1}
\int_N^{\infty} x^{-q/{p^*}- 1} \ll[ \frac{x^q}{(x^r+ d)^{q/r}} \rr] \d x \leq \int_1^{\infty} (x+ 1)^{-q/{p^*}- 1} \ll[ \frac{x^q}{(x^r+ d)^{q/r}} \rr] \d x. 
\de
In fact, we have
$$
\int_N^{\infty}  x^{-q/{p^*}- 1} \ll[ \frac{x^q}{(x^r+ d)^{q/r}} \rr] \d x \leq \int_N^{\infty} x^{-q/{p^*}- 1} \d x = \frac{p^*}{q} N^{-q/{p^*}}. 
$$
The existence of $N$ is obvious since the left side of (\ref{left1}) decreases to $0$ as $N \uparrow \infty$. Fix this sufficiently large integer $N$, then the left side of (\ref{left1}) is calculable. Using the integral transform $s^{-1}= d^{-1} x^r +1$, we have 
\be\lb{left2}
\int_N^{\infty} x^{-q/{p^*}- 1} \ll[ \frac{x^q}{(x^r+ d)^{q/r}} \rr] \d x = r^{-1} d^{- \frac{q}{rp^*}} B\ll(\frac{1+r}{r}, \frac{q- r- 1}{r}, \frac{d}{N^r+ d}\rr),
\de
where $B(a, b, x)$ is the incomplete Beta function: 
$$
B(a, b, x)= \int_{0}^{x} s^{a- 1} (1- s)^{b- 1} \d s.
$$
Applying the mean value theorem, (\ref{left1}) and (\ref{left2}), we have
\begin{align}\lb{left3}
\int_1^\infty & \ll[x^{-q/{p^*}}- (x+ 1)^{-q/{p^*}} \rr] \frac{x^q}{(x^r+ d)^{q/r}} \d x \notag \\
& \geq \int_1^\infty \frac{q}{p^*} (x+ 1)^{-q/{p^*}- 1} \ll[\frac{x^q}{(x^r+ d)^{q/r}} \rr] \d x \notag \\
& \geq d^{- \frac{q}{rp^*}} \ll(\frac{q}{p^*} \rr) r^{-1} B\ll(\frac{1+r}{r}, \frac{q- r- 1}{r}, \frac{d}{N^r+ d}\rr). 
\end{align}

Now, it's very easy to calculate the optimal constants. Using the relation
$$
B(a+ 1, b- 1)= \frac{a}{b- 1} B(a, b),
$$
it follows from (\ref{A}), (\ref{right}) and (\ref{left3}) that
$$
\aligned
A^q &\geq \ll[\sum_{n=1}^{\infty} \ll(n^{-q/p^*}- (n+ 1)^{-q/p^*} \rr) (H \mathbf{x} (n))^q \rr] \ll( \sum_{n=1}^\infty x_n^p \rr)^{- q/p} \\
&\geq \ll[\int_{1}^{\infty} \ll[x^{-q/p^*} - (x+ 1)^{-q/p^*} \rr] \frac{x^q}{(x^r+ d)^{\frac{q}{r}}} \d x \rr] \ll( \sum_{n=1}^\infty x_n^p \rr)^{- q/p} \\
&\geq \ll(\frac{q}{p^*}\rr) r^{q/p - 1} \cdot B\ll(\frac{1+ r}{r}, \frac{q- 1- r}{r}, \frac{d}{N^r+ d} \rr) \cdot B\ll(\frac{1}{r}, \frac{q- 1}{r}\rr)^{- q/p} \\
&\rightarrow k_{q,p}^q \qqd \text{(as $d \rightarrow \infty$)}. 
\endaligned
$$
Hence, the factor of basic upper estimate is best possible. 

(c) The final step is to remove condition $H \mathbf{\hat{v}} (\infty)= \infty$. We use Proposition \ref{C. estimates} of section \ref{Interval}. Fix $N_0 < \infty$. For given $\mathbf{u}$ and $\mathbf{v}$ on $[1, \infty)$, we define $\mathbf{u}^{N_0}$ and $\mathbf{v}^{N_0}$ to be the restriction of $\mathbf{u}$ and $\mathbf{v}$ on $[1, N_0]$. Then define 
$$
\overline{u}_n =
\begin{cases}
u_n^{N_0}, & 1\leq n \leq N_0, \\
0, & n> N_0, 
\end{cases}
$$
and 
$$
\overline{v}_n =
\begin{cases}
v_n^{N_0}, & 1\leq n\leq N_0, \\
1, & n> N_0.
\end{cases}
$$
Obviously, we have $H \mathbf{\overline{v}} (\infty) = \infty$. Applying the result of (a), we have 
$$
A(\mathbf{\overline{u}}, \mathbf{\overline{v}}) \leq k_{q, p} B(\mathbf{\overline{u}}, \mathbf{\overline{v}}).
$$
By Proposition \ref{C. estimates}, we get
$$
A(\mathbf{u}^{N_0}, \mathbf{v}^{N_0}) \leq k_{q, p} B(\mathbf{u}^{N_0}, \mathbf{v}^{N_0}).
$$
The assertion follows by letting $N_0\rightarrow \infty$. This completes the proof of Theorem \ref{kqpB} in the case $N= \infty$. 
\deprf

Review part (a) of the proof of Theorem (\ref{kqpB}), when 
\be\lb{v}
N= \infty, \qqd H \mathbf{\hat{v}} (\infty) =\infty,
\de
we give the method to construct $\mathbf{u}$ from $\mathbf{v}$ such that the Hardy-type inequalities (\ref{Hardy})  hold with these $\mathbf{u}$ and $\mathbf{v}$. The part (b) show the optimal constant reaches the upper bound of the basic estimate. It means that the basic upper estimate with the improved factor $k_{q,p}$ holds for a large class of $(\mathbf{u}, \mathbf{v})$. The original idea of this construction is from Chen \rf{Chen3}{Proposition 4.5}. To distinguish it from Theorem \ref{kqpB}, we give the following proposition. 

\prp\lb{Bliss_prp}
{\cms
For any positive sequences $\mathbf{v}$ and constant $0< C< \infty$, the discrete Hardy-type inequalities (\ref{Hardy}) hold on $[1, \infty)$ with 
\be
\tilde{u}_n = C^q \ll( (H \mathbf{\hat{v}} (n))^{-q/p^*}- (H \mathbf{\hat{v}} (n+1))^{-q/p^*} \rr), \qqd n\geq 1, 
\de
and its optimal constant $A$ satisfies 
\be
A \leq k_{q, p} C, 
\de
where $k_{q,p}$ is defined as (\ref{k_qp}). Moreover, when $N$ and $\mathbf{\hat{v}}$ satisfy (\ref{v}), the upper bound is sharp with $C= B$.
}
\deprp

\section{Hardy-type Inequalities on Interval}\lb{Interval}
In this section, we study the comparison results of the optimal constants and the basic estimates on different intervals. In continuous case, the corresponding comparison results have been done by Chen \rf{Chen3}{Appendix}. With these results, we can get the complete proofs of Theorem \ref{Var} and Theorem \ref{kqpB}. 

Before specific discussion, we need some notations. Fix two natural numbers $N$ and $N'$ with $N< N'$. Given two positive sequences $\mathbf{u}$ and $\mathbf{v}$ on $[1, N]$, we can extend them to $[1, N']$ as follows:
\be\lb{u'}
u'_i =
\begin{cases}
u_i, & 1\leq i\leq N, \\
0, & N< i\leq N';
\end{cases}
\de
\be\lb{v'}
v'_i =
\begin{cases}
v_i, & 1\leq i\leq N, \\
\#, & N< i\leq N',
\end{cases}
\de
where $\#$ means arbitrary positive numbers. Denote $A_N(\mathbf{u}, \mathbf{v})$ be the optimal constant of the Hardy-type inequalities (\ref{Hardy}) in the interval $[1, N]$ with sequences $\mathbf{u}$ and $\mathbf{v}$, and similar for $B_N(\mathbf{u}, \mathbf{v})$. 

The first result is a comparison for the optimal constants on different intervals. 

\prp\lb{C. constant}
{\cms
Given two positive sequences $\mathbf{u}'$ and $\mathbf{v}'$ on $[1, N']$. Use $\mathbf{u}$ and $\mathbf{v}$ to denote their restrictions to $[1, N]$. Then we have $A_N(\mathbf{u}, \mathbf{v}) \uparrow A_{N'}(\mathbf{u}', \mathbf{v}')$ as $N \uparrow {N'} \leq \infty$. 

In particular, if the inequality (\ref{Hardy}) holds on $[1, N']$, then it also holds with the same constant $A_{N'}(\mathbf{u}', \mathbf{v}')$ on $[1, N]$. 
}
\deprp

\prf
(a) Given a non-negative sequence $\mathbf{x}$ on $[1, N]$, we can extend to $[1, N']$ by setting 
\be\lb{extension of x}
x'_i =
\begin{cases}
x_i & 1\leq i\leq N, \\
0 & N< i\leq N'.
\end{cases} 
\de
Then we have 
$$
\aligned
\ll[ \sum_{n=1}^N u_n \ll(H \mathbf{x} (n) \rr)^q \rr]^{1/q} &= \ll[ \sum_{n=1}^{N'} u'_n \ll(H \mathbf{x'} (n) \rr)^q \rr]^{1/q} \\ 
&\leq A_{N'}(\mathbf{u}', \mathbf{v}') \ll[ \sum_{n=1}^{N'} v'_n {x'}_n^p \rr]^{1/p} \\
&= A_{N'}(\mathbf{u}', \mathbf{v}') \ll[ \sum_{n=1}^{N} v_n {x}_n^p \rr]^{1/p}.
\endaligned
$$
It means that $A_N(\mathbf{u}, \mathbf{v}) \leq A_{N'}(\mathbf{u}', \mathbf{v}')$. 

(b) Our next goal is to show the convergence. First we consider the case that $\sum_{n=1}^{N'} u'_n = \infty$. Clearly, in this case we have $N'= \infty$ and $A_{N'}(\mathbf{u}', \mathbf{v}')= \infty$. Besides, restricting to $[1, n]$ and choosing $\mathbf{x}= (1, 0, \dots, 0)$, we obtain 
$$
A_n(\mathbf{u}, \mathbf{v}) \geq \ll(\sum_{i=1}^n u_i \rr)^{1/q} v_1^{- 1/p} \rightarrow \infty, \qqd \text{as $n\rightarrow \infty$}.
$$
Hence the convergence holds in this case. 

(c) Let $\sum_{n=1}^{N'} u'_n < \infty$. For every non-negative sequence $\mathbf{x}$ on $[1, N']$ with $\sum_{n = 1}^{N'} v'_n x_n^p < \infty$, we get
$$
\frac{\ll[\sum_{n= 1}^N u_n (H \mathbf{x} (n))^q \rr]^{1/q}}{\ll(\sum_{n= 1}^N v_n x_n^p \rr)^{1/p}} \rightarrow \frac{\ll[\sum_{n= 1}^{N'} u'_n (H \mathbf{x} (n))^q \rr]^{1/q}}{ \ll( \sum_{n= 1}^{N'} v'_n x_n^p \rr)^{1/p}} \leq A_{N'}(\mathbf{u}', \mathbf{v}'),
$$
as $N \uparrow N'$. With (\ref{A}), for every $\vz> 0$, we can choose a sequence $\mathbf{x}$ such that
$$
A_{N'}(\mathbf{u}', \mathbf{v}') \leq \frac{\ll[\sum_{n= 1}^{N'} u'_n (H \mathbf{x} (n))^q \rr]^{1/q}}{\ll(\sum_{n= 1}^{N'} v'_n x_n^p \rr)^{1/p}} + \vz.
$$
Then we can choose $N$ closed to $N'$ such that
$$
\frac{\ll[\sum_{n= 1}^{N'} u'_n (H \mathbf{x} (n))^q \rr]^{1/q}}{\ll(\sum_{n= 1}^{N'} v'_n x_n^p \rr)^{1/p}} \leq \frac{\ll[\sum_{n= 1}^N u_n (H \mathbf{x} (n))^q \rr]^{1/q}}{\ll(\sum_{n= 1}^N v_n x_n^p \rr)^{1/p}} + \vz.
$$
Hence, we get
$$
A_N(\mathbf{u}, \mathbf{v}) \leq A_{N'}(\mathbf{u}', \mathbf{v}') \leq \frac{\ll[\sum_{n= 1}^N u_n (H \mathbf{x} (n))^q \rr]^{1/q}}{\ll(\sum_{n= 1}^N v_n x_n^p \rr)^{1/p}} + 2\vz \leq A_N(\mathbf{u}, \mathbf{v}) + 2\vz.
$$
It means that the convergence holds.
\deprf

The following result is about the factor in the basic estimates. 

\prp\lb{C. estimates}
{\cms
Given two positive sequences $\mathbf{u}$ and $\mathbf{v}$ on $[1, N]$, $\mathbf{u}'$ and $\mathbf{v}'$, defined by (\ref{u'}) and (\ref{v'}), are the extensions on $[1, N']$. Suppose that $A_{N'}(\mathbf{u}', \mathbf{v}')\leq k B_{N'}(\mathbf{u}', \mathbf{v}')$ for a universal constant $k$, then we have $A_N(\mathbf{u}, \mathbf{v})\leq k B_N(\mathbf{u}, \mathbf{v})$. 
}
\deprp

\prf
Given a sequence $\mathbf{x}$ on $[1, N]$, we can extend it from $[1, N]$ to $[1, N']$ by (\ref{extension of x}). Then we have
$$
\aligned
\ll( \sum_{n=1}^N u_n (H \mathbf{x} (n))^q \rr)^{1/q} &= \ll( \sum_{n=1}^{N'} u'_n (H \mathbf{x}' (n))^q \rr)^{1/q} \\
&\leq A_{N'}(\mathbf{u}', \mathbf{v}') \ll( \sum_{n= 1}^{N'} v'_n {x'}_n^p \rr)^{1/p} \\
&\leq k B_{N'}(\mathbf{u}', \mathbf{v}') \ll( \sum_{n= 1}^{N'} v'_n {x'}_n^p \rr)^{1/p} \\
&= k B_{N'}(\mathbf{u}', \mathbf{v}') \ll( \sum_{n= 1}^{N} v_n x_n^p \rr)^{1/p}. 
\endaligned
$$
With the definition of the extensions (\ref{u'}) and (\ref{v'}), we can easily check that 
$$
B_{N'}(\mathbf{u}', \mathbf{v}')= B_{N}(\mathbf{u}, \mathbf{v}).
$$
It follows that 
$$
\ll( \sum_{n=1}^N u_n (H \mathbf{x} (n))^q \rr)^{1/q} \leq k B_{N}(\mathbf{u}, \mathbf{v}) \ll( \sum_{n= 1}^{N} v_n x_n^p \rr)^{1/p}. 
$$
Hence $A_N(\mathbf{u}, \mathbf{v})\leq k B_N(\mathbf{u}, \mathbf{v})$ as required. 
\deprf

With the help of Proposition \ref{C. constant} and Proposition \ref{C. estimates}, we know the variational formulas of the optimal constants, the basic estimates and the improved factor of the basic upper estimates are true when $N< \infty$. So far, we complete the proofs of our main results. 

The following result gives an opposite view of Proposition \ref{C. constant}: from some local sub-intervals to the whole interval. It gives us an approximating procedure for the unbounded interval. 

\prp\lb{C. interval}
{\cms
Given two positive sequences $\mathbf{u}$ and $\mathbf{v}$ on $[1, N]$, extend them to $[1, N']$ by (\ref{u'}) and (\ref{v'}). Then we have $A_N (\mathbf{u}, \mathbf{v}) = A_{N'} (\mathbf{u}', \mathbf{v}')$. 
}
\deprp

\prf
For any sequence $\mathbf{x} \in \aq[1, N]$, let $\mathbf{x}'$ be the extension of $\mathbf{x}$ from $[1, N]$ to $[1, N']$ by (\ref{extension of x}). Obviously, we have $\mathbf{x}' \in \aq[1, N']$. The inequalities in $[1, N']$ are
$$
\| H \mathbf{x}' \|_{l^q(u')} \leq A_{N'} (\mathbf{u}', \mathbf{v}') \| \mathbf{x}' \|_{l^p(v')}.
$$
With (\ref{u'}), (\ref{v'}) and (\ref{extension of x}), it follows that
$$
\| H \mathbf{x} \|_{l^q(u)} \leq A_{N'} (\mathbf{u}', \mathbf{v}') \| \mathbf{x} \|_{l^p(v)}.  
$$
Because $\mathbf{x}$ is arbitrary, it implies that $A_N (\mathbf{u}, \mathbf{v}) \leq A_{N'} (\mathbf{u}', \mathbf{v}')$. 

Conversely, for any $\mathbf{x} \in \aq[1, N']$, we have
$$
\aligned
\ll(\sum_{n=1}^{N'} {u'}_n (H \mathbf{x} (n))^q \rr)^{1/q} &= \ll(\sum_{n=1}^{N} u_n (H \mathbf{x} (n))^q \rr)^{1/q} \\
&\leq A_N (\mathbf{u}, \mathbf{v}) \ll( \sum_{n=1}^N v_n \mathbf{x}_n^p \rr)^{1/p} \\
&\leq A_N (\mathbf{u}, \mathbf{v}) \ll( \sum_{n=1}^{N'} v'_n \mathbf{x}_n^p \rr)^{1/p}. 
\endaligned
$$
This implies that $A_{N'} (\mathbf{u}', \mathbf{v}') \leq A_N (\mathbf{u}, \mathbf{v})$ and then the equality holds.
\deprf

\section{Examples}\lb{examples}

As mentioned in introduction, Hardy-type inequalities play important role in probability theory. The first example is from birth-death processes which is standard having constant birth and death rates, cf. \rf{Chen5}{Example 5.3}. We present this example to illustrate the power of out results. 
\xmp\lb{example1}
{\cms
Let $p=q=2$ and $N= \infty$. For $n \geq 1$, let $u_n = \gz^n$, $v_n= b \gz^n$, where $\gz$ and $b$ are constants with $\gz< 1$ and $b> 0$. Then
$$
B < \widetilde{\dz}_1= \overline{\dz}_1 < A = \dz_1 < 2B,
$$
where $B= \dfrac{1}{\sqrt{b} (1- \gz)}$, $\widetilde{\dz}_1= \overline{\dz}_1= \dfrac{\sqrt{1+ \gz}}{\sqrt{b} (1- \gz)}$, $A = \dz_1= \dfrac{1}{\sqrt{b}(1- \sqrt{\gz})}$. Moreover, the optimal constant is attained at sequence
$$
a_n= \gz^{(-n +1)/2} \ll[n - (n- 1) \gz^{1/2} \rr], \qqd n\geq 1.
$$
}
\dexmp

\prf (a) First, $B$ is easy to calculate. By the definition, we have
$$
B= \sup_{n\in [1, \infty)}  \ll(\sum_{i=1}^n b^{-1} \gz^{- i} \rr)^{1/2} \ll(\sum_{j= n}^{\infty} \gz^{j}\rr)^{1/2} = \frac{1}{\sqrt{b} (1- \gz)}. 
$$
Next, by (\ref{k_pp}), we have $k_{2, 2}= 2$. By Corollary \ref{basic estimates}, we obtain the basic estimates of the optimal constants: 
\be\lb{EX13}
\frac{1}{\sqrt{b} (1- \gz)} \leq A \leq \frac{2}{\sqrt{b} (1- \gz)}. 
\de

(b) To compute $\dz_1$, we use Corollary \ref{approximating}. Let 
$$
x_n^{(1)}= (H \mathbf{\hat{v}} (n))^{1/2}- (H \mathbf{\hat{v}} (n- 1))^{1/2}, \qqd n\geq 1,
$$
and then
$$
H \mathbf{x}^{(1)} (n) = (H \mathbf{\hat{v}} (n))^{1/2} = \ll[\frac{\gz^{-n}- 1}{b (1- \gz)}\rr]^{1/2}. 
$$
For convenience, we use $\fz_n= \gz^{-n} -1$ in the following. By direct computations, we have
\begin{align}\lb{EX11}
I\!I_n^* \ll(\mathbf{x}^{(1)}\rr) &= \frac{1}{H \mathbf{x}^{(1)} (n)} \sum_{i=1}^n \hat{v}_i \ll(\sum_{j=i}^\infty u_j \ll(H \mathbf{x}^{(1)} (j) \rr) \rr) \notag \\
&= \frac{b^{-3/2}}{(1- \gz)^{1/2}} \frac{1}{H \mathbf{x}^{(1)} (n)} \sum_{i=1}^n \gz^{-i}  \ll(\sum_{j=i}^\infty \gz^j \fz_j^{1/2} \rr) \notag \\
&= \frac{\fz_n^{-1/2}}{b (1- \gz)} \ll[\sum_{j=1}^n \gz^{j} \fz_j^{3/2} + \fz_n \sum_{j=n+ 1}^\infty \gz^{j} \fz_j^{1/2}\rr]. 
\end{align}
At the last step, we exchange the order of summation. From (\ref{EX11}), it is easy to check that $I\!I_n^* \ll(\mathbf{x}^{(1)}\rr)$ reaches the maximum when $n \rightarrow \infty$. Hence, by L'Hospital's rule, we obtain
$$
\aligned
\dz_{1}^2 &= \sup_{n\in [1, \infty)} I\!I_n^* \ll(\mathbf{x}^{(1)} \rr) \\
&= \frac{1}{b (1- \gz)} \ll[ \lim_{n\rightarrow \infty} \fz_n^{-1/2} \sum_{j=1}^n \gz^{j} \fz_j^{3/2} + \lim_{n\rightarrow \infty} \fz_n^{1/2} \sum_{j=n+ 1}^\infty \gz^{j} \fz_j^{1/2} \rr] \\
&= \frac{1}{b (1- \gz)} \ll[ \frac{1}{1- \sqrt{\gz} } + \frac{\sqrt{\gz}}{1- \sqrt{\gz}} \rr] \\
&= \frac{1}{b (1- \sqrt{\gz})^2}. 
\endaligned
$$

(c) Similarly, we use Corollary \ref{approximating} to compute $\overline{\dz}_1$ and $\widetilde{\dz}_1$. Fix $k> 0$, let 
$$
y_n^{(k, 1)}= 
\begin{cases}
b^{-1} \gz^{-n}, & n\leq k, \\
0, & n>k,
\end{cases}
$$
and then 
$$
H \mathbf{y}^{(k, 1)} (n) = \frac{\gz^{-(n \wedge k)} - 1}{b (1- \gz)}= \frac{\fz_{n \wedge k}}{b (1- \gz)}.
$$
By lots of tedious calculations, we have 
$$
\aligned
I\!I_n \ll( \mathbf{y}^{(k, 1)} \rr) &= \frac{1}{H \mathbf{y}^{(k, 1)} (n)} \sum_{i=1}^n \hat{v}_i \ll(\sum_{j= i}^\infty u_j \ll(H \mathbf{y}^{(k, 1)} \rr) \rr) \\
&= \frac{1}{b \fz_{n \wedge k}} \sum_{i=1}^n \gz^{-i} \ll(\sum_{j= i}^\infty \gz^{j} \fz_{j \wedge k} \rr) \\
&= \frac{1}{b (1- \gz)} \ll[\frac{1+ \gz}{1- \gz} - \frac{2 (n \wedge k)}{\fz_{n \wedge k}} + (k- n)\vee 0 - \frac{\gz^{k+ 1} \fz_{(n- k) \vee 0}}{1- \gz}  \rr].
\endaligned
$$
Next, note that $I\!I_n \ll( \mathbf{y}^{(k, 1)} \rr)$ reaches the minimum when $n= k$, and then
$$
\aligned
\widetilde{\dz}_1^2 &= \sup_{k\in [1, \infty)} \inf_{n\in [1, \infty)} I\!I_n \ll(\mathbf{y}^{(k, 1)} \rr) \\
&= \sup_{k\in [1, \infty)} \frac{1}{b (1- \gz)} \ll(\frac{1+ \gz}{1- \gz} - \frac{2 k}{\fz_{k}} \rr) \\
&= \frac{1}{b (1- \gz)} \lim_{k\rightarrow \infty} \ll(\frac{1+ \gz}{1- \gz} - \frac{2 k}{\fz_{k}} \rr) \\
&= \frac{1+ \gz}{b (1- \gz)^2}.
\endaligned
$$

Now, we consider $\overline{\dz}_1$. Since
$$
\sum_{n=1}^\infty v_n \ll(y_n^{(k, 1)} \rr)^2 = \frac{\fz_k}{b(1- \gz)}, 
$$
and 
$$
\sum_{n=1}^\infty u_n \ll( H \mathbf{y}^{(k, 1)} (n) \rr)^2 = \frac{1}{b^2 (1- \gz)^2} \ll[\sum_{n= 1}^k \gz^n \fz_n^2 + \frac{\gz^{k+ 1} \fz_k^2}{1- \gz} \rr], 
$$
we have 
$$
\aligned
\overline{\dz}_1^2&= \sup_{k\in [1, \infty)} \frac{1}{b(1- \gz)} \ll[\fz_k^{-1} \sum_{n=1}^k \gz^n \fz_n^2 + \frac{\gz- \gz^{k+ 1}}{1- \gz}  \rr] \\
&= \frac{1}{b(1- \gz)} \ll[ \lim_{k\rightarrow \infty} \fz_k^{-1} \sum_{n=1}^k \gz^n \fz_n^2 + \frac{\gz }{1- \gz}  \rr] = \frac{1+ \gz}{b (1- \gz)^2}.
\endaligned
$$
In the last step, the L'Hospital's rule is used to calculate the limitation of $k$. 

(d) So far, by Corollary \ref{approximating}, we obtain the estimates of the optimal constants,  which is more precise than the basic estimates (\ref{EX13}) 
\be\lb{EX12}
\frac{\sqrt{1+ \gz}}{\sqrt{b} (1- \gz)} \leq A \leq \frac{1}{\sqrt{b} (1- \sqrt{\gz})}. 
\de

In fact, the optimal constant can be accurately calculated. Let $a_n= \gz^{(-n +1)/2} \ll[n - (n- 1) \gz^{1/2} \rr] (n\geq 1)$, then
$$
H \mathbf{a} (n)= n \gz^{(-n+ 1)/2}. 
$$
Here we want to use $\mathbf{a}$ instead of $\mathbf{y}^{(k, 1)}$ to get the lower estimates. However, it is easy to check that $\mathbf{a}$ is non-summability. It means that Theorem \ref{Var} is invalid. By the classical variational formula (\ref{A}) and the L'Hospital's rule, we have
$$
\aligned
A^2 &\geq \frac{\sum_{n= 1}^\infty \gz^n a_n^2 }{\sum_{n=1}^\infty b \gz^n \ll(a_n - a_{n-1} \rr)^2} \\
&= b^{-1} \lim_{n\rightarrow \infty} \frac{n^2}{\ll[n- (n- 1) \gz^{1/2} \rr]^2} \\
&= \frac{1}{b (1- \sqrt{\gz})^2}.
\endaligned
$$
As a consequence, we obtain $A = \dfrac{1}{\sqrt{b} (1- \sqrt{\gz})}$.
\deprf

To distinguish the first example, the second one is about the nonlinear situation $p \neq q$, which is from proof of Theorem \ref{kqpB}. The optimal constant is clear in this example.  
\xmp\lb{example2}
{\cms
Let $p \neq q$ and $N= \infty$. For $n\geq 1$, let $u_n= n^{-q/{p^*}}- (n+ 1)^{-q/{p^*}}$, $v_n \equiv 1$. Then

(1) The optimal constant is $A= k_{q, p}$, which is attained at sequence $\mathbf{x}$: 
$$
x_n= \frac{cn}{(n^r+ d)^{1/r}} - \frac{c(n- 1)}{((n-1)^r+ d)^{1/r}}, \qqd n \geq 1,
$$
where $r=q/p - 1$, $k_{q, p}$ is defined as (\ref{k_qp}), $c$ and $d$ are arbitrary positive constants. 

(2) The basic estimates and the approximating procedure are
$$
B \leq \overline{\dz}_1 \vee \widetilde{\dz}_1 \leq A = k_{q, p} B \leq \dz_1, 
$$ 
where $B=1$, $\overline{\dz}_1 \geq 1$, $\widetilde{\dz}_1 \geq 1$ and $\dz_1 \leq \ll(1+ \dfrac{q}{p^*} \rr)^{1/q+ 1/p^*}$.
}
\dexmp

\prf The first part has been done in Theorem \ref{kqpB}. The remainder of this proof is to compute $\dz_1$, $\overline{\dz}_1$ and $\widetilde{\dz}_1$. 

To compute $\dz_1$, let
$$
x_n^{(1)}= n^{q/(p^*+ q)} - (n- 1)^{q/(p^*+ q)}, n\geq 1, 
$$
then we have 
$$
\aligned
I\!I_n^* \ll(\mathbf{x}^{(1)}\rr) &= n^{- \frac{q}{p^*+q}} \sum_{i= 1}^n \ll\{\sum_{j=i}^{\infty} \ll[j^{-\frac{q}{p^*}}- (j+ 1)^{- \frac{q}{p^*}} \rr] j^{\frac{q^2}{p^*(p^*+ q)}} \rr\}^{p^*/q} \\
&\leq n^{- \frac{q}{p^*+q}} \sum_{i= 1}^n \ll\{ \ll(\frac{q}{p^*} \rr) \sum_{j=i}^{\infty} \int_j^{j+1} x^{\frac{q^2}{p^*(p^*+ q)}- \frac{q}{p^*}- 1} \d x \rr\}^{p^*/q} \\
&= n^{- \frac{q}{p^*+q}} \sum_{i= 1}^n \ll[ \ll(\frac{p^*+ q}{p^*} \rr) i^{- \frac{q}{p^*+ q}} \rr]^{p^*/q} \\
&\leq n^{- \frac{q}{p^*+q}} \ll(\frac{p^*+ q}{p^*} \rr)^{p^*/q} \ll(1+ \int_1^n x^{- \frac{p^*}{p^*+ q}} \d x\rr) \\
&= \ll(1+ \frac{q}{p^*} \rr)^{p^*/q + 1}. \\
\endaligned
$$
Therefore, we obtain
$$
\dz_1 = \sup_{n\in [1, \infty)} \ll[ I\!I_n^* \ll(\mathbf{x}^{(1)} \rr)\rr]^{1/p^*} \leq \ll(1+ \frac{q}{p^*} \rr)^{1/q + 1/p^*}.
$$

To compute $\overline{\dz}_1$ and $\widetilde{\dz}_1$, let
$$
y_n^{(k, 1)} =
\begin{cases}
1, & n\leq k, \\
0, & n> k. 
\end{cases}
$$
Obviously, we have $H \mathbf{y}^{(k, 1)} (n) = n \wedge k$, $\|\mathbf{y}^{(k, 1)} \|_{l^p(v)} = k^{1/p}$ and 
$$
\aligned
\|H \mathbf{y}^{(k, 1)} \|_{l^q(u)} &= \ll[ \sum_{n=1}^\infty u_n (n \wedge k)^q \rr]^{1/q} \\
&= \ll[ \sum_{n=1}^{k- 1} \ll(n^{-\frac{q}{p^*}}- (n+1)^{-\frac{q}{p^*}} \rr) n^q + k^{q/p} \rr]^{1/q}.
\endaligned
$$
Hence, we obtain
$$
\aligned
\overline{\dz}_1 &= \sup_{k\in [1, \infty)} \frac{\|H \mathbf{y}^{(k, 1)} \|_{l^q(u)}}{\|\mathbf{y}^{(k, 1)} \|_{l^p(v)}} \\
&= \sup_{k\in [1, \infty)} \ll[k^{-q/p} \sum_{n=1}^{k- 1} \ll(n^{-q/p^*}- (n+1)^{-q/p^*} \rr) n^q + 1 \rr]^{1/q} \\ 
&\geq 1.
\endaligned
$$
Now, we consider $\widetilde{\dz}_1$. With directly calculating, we have
$$
\aligned
I\!I_n & \ll(\mathbf{y}^{(k, 1)}\rr) = \inf_{n\in [1, \infty)} \frac{1}{n\wedge k} \sum_{i=1}^n \ll[\sum_{j=i}^{\infty} u_j \ll(j \wedge k\rr)^{q- 1} \rr]^{p^*- 1}\\
&= \frac{1}{n \wedge k} \sum_{i=1}^{k \wedge n} \ll[k^{q/p - 1} + \sum_{j=i}^{k-1} u_j j^{q-1} \rr]^{p^*- 1} + \mathbbm{1}_{\{n> k\}} k^{\frac{q- p}{p- 1}} \sum_{i=k +1}^n i^{-q/p}. 
\endaligned
$$
Obviously, $I\!I_n \ll(\mathbf{y}^{(k, 1)}\rr)$ is increasing when $n\geq k$. It means that $I\!I_n \ll(\mathbf{y}^{(k, 1)}\rr)$ reaches its minimum at $n\in [1, k]$. Thus, we obtain
$$
\aligned
\inf_{n\in [1, \infty)} I\!I_n \ll(\mathbf{y}^{(k, 1)}\rr) &= \inf_{n\leq k} \frac{1}{n} \sum_{i=1}^{n} \ll[k^{q/p - 1} + \sum_{j=i}^{k-1} u_j j^{q-1} \rr]^{p^*- 1} \\
&\geq \inf_{n\leq k} \ll[k^{q/p - 1} + \sum_{j=n}^{k-1} u_j j^{q-1} \rr]^{p^*- 1} \\
&= k^{(q/p - 1)(p^*- 1)}.
\endaligned
$$
Therefore, we obtain
$$
\aligned
\widetilde{\dz}_1 &= \sup_{k\in [1, \infty)} k^{1/q- 1/p} \ll( \inf_{n\in [1, \infty)} I\!I_n \ll(\mathbf{y}^{(k, 1)}\rr) \rr)^{(p- 1)/q} \\
&\geq \sup_{k\in [1, \infty)} k^{1/q- 1/p} \ll[k^{(q/p - 1)(p^*- 1)} \rr]^{(p- 1)/q}= 1. \qqd \square
\endaligned
$$

\noindent {\bf Acknowledgements} $\qd$ This paper is based on the series of studies of my supervisor Prof. M. F. Chen. Heartfelt thanks are given to my supervisor for his careful guidance and helpful suggestions. Thanks are also given to Prof. Y. H. Mao, Prof. F. Y. Wang and Prof. Y. H. Zhang for their comments and suggestions, which lead to lots of improvements of this paper. 

The research is supported by NSFC (Grant No. 11131003) and by the ``985'' project from the Ministry of Education in China.

\end{document}